\documentclass[12pt]{article}
\usepackage{amsmath,amssymb,url}
\usepackage[none]{hyphenat}
\begin{document}

\begin{center}{\Large \bf Five families of rapidly convergent evaluations of zeta values}\\[3pt]
{\large David Broadhurst, Open University, 3 April 2024}\end{center}{\large

\noindent This work derives 5 methods to evaluate families of odd zeta values by
combining a power of $\pi$ with Lambert series whose ratios of successive terms
tend to $e^{-\pi\sqrt{a}}$ with integers $a\ge7$, outperforming Ramanujan's results
with merit $a=4$. Families with $a=7$ and $a=8$ evaluate $\zeta(2n+1)$.
Families with $a=9$ and $a=16$ evaluate $\zeta(4n+1)$ with faster convergence.
A fifth family with $a=12$ evaluates $\zeta(6n+1)$ and gives the fastest convergence for $\zeta(6n+7)$.
Members of three of the families were discovered empirically by Simon Plouffe.
An intensive new search strongly suggests that there are no more than 5 families with integers $a\ge7$.
There are at least 20 families that involve Lambert series with rational $a>4$.
Quasi-modular transformations of Lambert series resolve rational
sequences that were discovered empirically. Expansions of Lambert series in
polylogarithms, familiar from quantum field theory, provide
proofs of all known evaluations with rational merit $a>4$.

\section{Introduction}

I derive 5 families of reductions of $\zeta(k)$, for odd $k>1$,
to algebraic multiples of $\pi^k$ and rational multiples of very rapidly convergent Lambert series of the form
\begin{equation}S_k(q)=\sum_{n=1}^\infty\frac{1}{n^k}
\left(\frac{q^n}{1-q^n}\right)=\sum_{m=1}^\infty{\rm Li}_k(q^m)\label{Lam}\end{equation}
with expansions in polylogarithms ${\rm Li}_k(x)=\sum_{n>0}x^nn^{-k}$~\cite{Lewin}.
For each of these 5 families, the Lambert series have spectacularly fast convergence,
with arguments $q=e^{-\pi\sqrt{a}}<\frac{1}{4071}$ for integers $a\ge7$. All 5 families
give evaluations that converge faster than the classical evaluations~\cite{Bruce,Henri,P,V}
\begin{align}
\zeta(3)&=\frac{7\pi^3}{180}
-2\sum_{n=1}^\infty\frac{1}{n^3}\left(\frac{1}{e^{2\pi n}-1}\right)\label{rz3}\\
\zeta(5)&=\frac{\pi^5}{294}
-\frac{2}{35}\sum_{n=1}^\infty\frac{1}{n^5}\left(\frac{36}{e^{2\pi n}-1}
+\frac{1}{e^{2\pi n}+1}\right)\label{rz5}
\end{align}
with lesser merit $a=4$ and hence $q=e^{-2\pi}\approx\frac{1}{535}$ in Ramanujan's formula~(\ref{rz3}).

Section~2 describes the families. For each, I determine the rational
coefficients of Lambert series, by solving recurrence relations for integer sequences.
Section~3 comments on computational efficiency. Section~4 gives proofs of
the results stated in Section~2 and also derives the coefficients of powers of $\pi$.
Section~5 offers comments and conclusions.

\section{Families that outperform Ramanujan}

I assign merit $a$ to a series with an asymptotic ratio of terms $e^{-\pi\sqrt{a}}$.
Evaluations~(\ref{rz3},\ref{rz5}) have merit $a=4$. There are 5 families that have integer merit $a\ge7$,
with $a=7,8,9,16,12$. The respective square-free divisors of $a$ are $s=7,2,1,1,3$
and the Lambert series are of the form $S_k(q_s^d)$ where $q_s=e^{-\pi\sqrt{s}}$ and $d|m$ is a divisor of a least
common multiple $m=4,6,12,20,12$, with $d\sqrt{s}>2$ to ensure convergence faster
than Ramanujan's result~(\ref{rz3}).
I conjecture that there are precisely 5 families with integers $a\ge7$ and exhibit
15 families with Lambert series of fractional merit $a>4$. The order of merit
of the 20 families with rational $a>4$ is given in~(\ref{om}).

\subsection{Family A: evaluation of $\zeta(2n+1)$ with merit 7 using divisors of 4}

For odd $k>1$, there is an infinite family of evaluations of the form
\begin{equation}
\zeta(k)=A_{k,0}\sqrt7\pi^k+\sum_{d\in\{1,2,4\}}A_{k,d}S_k(q_7^d),\quad q_7=\exp(-\pi\sqrt7).\label{Ak}\end{equation}
Using the procedure {\tt lindep} in {\tt Pari/GP}~\cite{GP}, I obtained
\begin{align}
\zeta(3)&=\tfrac{29\sqrt7}{1980}\pi^3+\tfrac{24}{11}S_3(q_7)-\tfrac{52}{11}S_3(q_7^2)+\tfrac{6}{11}S_3(q_7^4)\label{A3}\\
\zeta(5)&=\tfrac{5\sqrt7}{3906}\pi^5+\tfrac{64}{31}S_5(q_7)-\tfrac{130}{31}S_5(q_7^2)+\tfrac{4}{31}S_5(q_7^4)\label{A5}\\
\zeta(7)&=\tfrac{851\sqrt7}{6747300}\pi^7+\tfrac{240}{119}S_7(q_7)-\tfrac{1927}{476}S_7(q_7^2)+\tfrac{15}{476}S_7(q_7^4)\label{A7}
\end{align}
and extended these results to $\zeta(127)$.

The existence of this family was inferred by Simon Plouffe~\cite{P1}. I was then able
to determine that the coefficients of the Lambert series in~(\ref{Ak}) satisfy
\begin{equation}A_{2n+1,1}=2^{2n} A_{2n+1,4}=2+\frac{2}{a_n}\,,\quad \sum_{d\in\{1,2,4\}}A_{2n+1,d}=-2\label{Ad}\end{equation}
with $a_n={\tt 11,\, 31,\, 119,\, 543,\, 1991,\, 8239,\, 32855,\, 130623,\, 525287,\, 2095951,\ldots}$\\
This integer sequence obeys the recurrence relation
\begin{equation}a_{n+3}=a_{n+2}+8a_{n+1}+16a_n\label{ar}\end{equation}
for $n\ge0$, with $a_0=0$. Observing that the
roots of $x^3=(x+4)^2$ are $x=4$ and $x=\frac12(-3\pm\sqrt{-7})$, I inferred that
\begin{equation}a_n=2^{2n+1}-\left(\frac{1+\sqrt{-7}}{2}\right)^{2n}-\left(\frac{1-\sqrt{-7}}{2}\right)^{2n}.\label{an}\end{equation}
Using the coefficient of $\sqrt7\pi^k$ in~(\ref{Ak}), I define the rational sequence
\begin{equation}A_n=\frac{(2n+3)!}{2^{2n-1}}a_n A_{2n+1,0}\label{Ai}\end{equation}
with $A_n={\tt 29/3,\, 25,\, 851/5,\, 6451/3,\, 1088813/35,\, 684521,\, 846243643/45,\ldots}$\\
{\bf Conjecture 1:} The denominator of $A_n$ is divisible by no prime greater than $n+2$.\\
{\bf Comment 1:} For $\zeta(291)$, I found that $3\cdot5\cdot13\cdot19\cdot59A_{145}$ is a 459-digit prime.
For $\zeta(2821)$, the numerator of $A_{1410}$ is divisible by a 7128-digit probable prime.

\subsection{Family B: evaluation of $\zeta(2n+1)$ with merit 8 using divisors of 6}

Following work by Linas Vep\v{s}tas~\cite{V}, Plouffe found faster convergence in evaluations of the form~\cite{P1}
\begin{equation}
\zeta(k)=B_{k,0}\sqrt2\pi^k+\sum_{d\in\{2,3,6\}}B_{k,d}S_k(q_2^d),\quad q_2=\exp(-\pi\sqrt2).\label{Bk}\end{equation}
Then {\tt lindep} gave
\begin{align}
\zeta(3)&=\tfrac{17\sqrt2}{620}\pi^3-\tfrac{60}{31}S_3(q_2^2)-\tfrac{4}{31}S_3(q_2^3)+\tfrac{2}{31}S_3(q_2^6)\label{B3}\\
\zeta(5)&=\tfrac{191\sqrt2}{79695}\pi^5-\tfrac{516}{253}S_5(q_2^2)+\tfrac{8}{253}S_5(q_2^3)+\tfrac{2}{253}S_5(q_2^6)\label{B5}\\
\zeta(7)&=\tfrac{3197\sqrt2}{13538700}\pi^7-\tfrac{612}{307}S_7(q_2^2)-\tfrac{16}{2149}S_7(q_2^3)+\tfrac{2}{2149}S_7(q_2^6)\label{B7}
\end{align}
and results up to $\zeta(127)$.

The coefficients of the Lambert series in~(\ref{Bk}) satisfy
\begin{equation}B_{2n+1,6}=\frac{B_{2n+1,3}}{(-2)^n}=\frac{2}{b_n}\,,\quad \sum_{d\in\{2,3,6\}}B_{2n+1,d}=-2\label{Bd}\end{equation}
with $b_n={\tt 31,\, 253,\, 2149,\, 19633,\, 177661,\, 1593601,\, 14346013,\, 129151873,\ldots}$\\
This integer sequence obeys the recurrence relation
\begin{equation}b_{n+4}=5b_{n+3}+23b_{n+2}+99b_{n+1}+162b_n\label{br}\end{equation}
for $n\ge0$, with $b_0=0$. Observing that the roots of $x^4=5x^3+23x^2+99x+162$
are $x=9$, $x=-2$ and $x=-1\pm\sqrt{-8}$, I inferred that
\begin{equation}b_n=3^{2n+1}-(-2)^n-(1+\sqrt{-2})^{2n}-(1-\sqrt{-2})^{2n}.\label{bn}\end{equation}
Using the coefficient of $\sqrt2\pi^k$ in~(\ref{Bk}), I define the rational sequence
\begin{equation}B_n=\frac{(2n+3)!}{2^{2n-1}}b_n B_{2n+1,0}\label{Bi}\end{equation}
with $B_n={\tt 51,\, 382,\, 28773/5,\, 145566,\, 181903711/35,\, 247890198,\ldots}$\\
{\bf Conjecture 2:} The denominator of $B_n$ is divisible by no prime greater than $n+2$.\\
{\bf Comment 2:} For $\zeta(189)$, I found that $77B_{94}/3298100598$ is a 287-digit prime.
For $\zeta(619)$, the numerator of $B_{309}$ is divisible by a 1268-digit prime.

\subsection{Family C: evaluation of $\zeta(4n+1)$ with merit 9 using divisors of 12}

For $k = 4n+1$, I found faster convergence for evaluations of the form
\begin{equation}
\zeta(k)=C_{k,0}\pi^k+\sum_{d\in\{3,4,6,12\}}C_{k,d}S_k(q_1^d),\quad q_1=\exp(-\pi).\label{Ck}\end{equation}
Then {\tt lindep} gave
\begin{align}
\zeta(5)&=\tfrac{682}{201285}\pi^5+\tfrac{296}{355}S_5(q_1^3)-\tfrac{488}{355}S_5(q_1^4)-\tfrac{1073}{710}S_5(q_1^6)
+\tfrac{37}{710}S_5(q_1^{12})\label{C5}\\
\zeta(9)&=\tfrac{5048}{150155775}\pi^9-\tfrac{2272}{1605}S_9(q_1^3)-\tfrac{5624}{1605}S_9(q_1^4)+\tfrac{37559}{12840}S_9(q_1^6)
-\tfrac{71}{12840}S_9(q_1^{12})\label{C9}
\end{align}
and results up to $\zeta(125)$.

The coefficients of the Lambert series in~(\ref{Ck}) satisfy
\begin{gather}C_{4n+1,12}=\frac{C_{4n+1,3}}{2^{4n}}=\left(1-\frac{2^{4n+1}+1}{(-4)^n}\right)\frac{2}{c_n},\label{Cd}\\
\frac{C_{4n+1,6}}{C_{4n+1,12}}=-(2^{4n+1}+(-4)^n+1),\quad\sum_{d\in\{3,4,6,12\}}C_{4n+1,d}=-2\label{Cd1}\end{gather}
with $c_n={\tt 355,\, 11235,\, 2114515,\, 95520195,\, 12606788275,\, 709832878755,\ldots}$\\
This integer sequence obeys the recurrence relation
\begin{equation}c_{n+3}=33c_{n+2}+4192c_{n+1}-82944c_n\label{cr}\end{equation}
for $n\ge0$, with $5|c_n$ and $c_0=0$. Observing that the roots of $x^3=33x^2+4912x-82944$
are $x=81$, $x=-64$ and $x=16$, I inferred that
\begin{equation}c_n=3^{4n+1}-2(-4)^{3n}-2^{4n}.\label{cn}\end{equation}
Using the coefficient of $\pi^k$ in~(\ref{Ck}), I define the rational sequence
\begin{equation}C_n=\frac{(4n)!(2n+1)(4n+3)}{2^{4n}}c_nC_{4n+1,0}\label{Ci}\end{equation}
with $C_n={\tt 341/9,\, 88340/27,\, 26827985/3,\, 18438055674,\ldots}$\\
{\bf Conjecture 3:} The denominator of $C_n$ is divisible by no prime greater than $2n+1$
and by no prime congruent to 1 modulo 4.\\
{\bf Comment 3:} For $\zeta(285)$, I found that $19\cdot23\cdot27C_{71}/28362151103$ is a 483-digit prime.
For $\zeta(1449)$, the numerator of $C_{362}$ is divisible by a 3495-digit probable prime.

\subsection{Family D: evaluation of $\zeta(4n+1)$ with merit 16 using divisors of 20}

For $k = 4n+1$, Plouffe found even faster convergence for evaluations of the form~\cite{P1}
\begin{equation}
\zeta(k)=D_{k,0}\pi^k+\sum_{d\in\{4,5,10,20\}}D_{k,d}S_k(q_1^d),\quad q_1=\exp(-\pi).\label{Dk}\end{equation}
Then {\tt lindep} gave
\begin{equation}
\zeta(5)=\tfrac{694}{204813}\pi^5-\tfrac{6280}{3251}S_5(q_1^4)+\tfrac{296}{3251}S_5(q_1^5)-\tfrac{1073}{6502}S_5(q_1^{10})
+\tfrac{37}{6502}S_5(q_1^{20})\label{D5}
\end{equation}
and results up to $\zeta(125)$.

The coefficients of the Lambert series in~(\ref{Dk}) satisfy
\begin{gather}D_{4n+1,20}=\frac{D_{4n+1,5}}{2^{4n}}=\left(1-\frac{2^{4n+1}+1}{(-4)^n}\right)\frac{2}{d_n},\label{Dd}\\
\frac{D_{4n+1,10}}{D_{4n+1,20}}=-(2^{4n+1}+(-4)^n+1),\quad\sum_{d\in\{4,5,10,20\}}D_{4n+1,d}=-2\label{Dd1}\end{gather}
with $d_n={\tt 3251,\, 1945731,\, 1221199811,\, 762905503491,\, 476839323944771,\ldots}$\\
This integer sequence obeys a fifth order recurrence relation whose characteristic polynomial is
$(x-625)(x+64)(x-16)(x^2+14x+625)$, with $d_0=0$, from which I inferred that
\begin{equation}d_n=5^{4n+1}-2(-4)^{3n}-2^{4n}-(2+\sqrt{-1})^{4n}-(2-\sqrt{-1})^{4n}.\label{dn}\end{equation}
Using the coefficient of $\pi^k$ in~(\ref{Dk}), I define the rational sequence
\begin{equation}D_n=\frac{(4n)!(2n+1)(4n+3)}{2^{4n}}d_n D_{4n+1,0}\label{Di}\end{equation}
with $D_n={\tt 347,\, 15297320/27,\, 5164889285,\, 1030773088821534/7,\ldots}$\\
{\bf Conjecture 4:} The denominator of $D_n$ is divisible by no prime greater than $2n-1$
and by no prime congruent to 1 modulo 4.\\
{\bf Comment 4:} For $\zeta(273)$, I found that $19\cdot23\cdot31\cdot21^2D_{68}/6706150$
is a 524-digit prime.
For $\zeta(897)$, the numerator of $D_{224}$ is divisible by a 2173-digit probable prime.

\subsection{Family E: evaluation of $\zeta(6n+1)$ with merit 12 using divisors of 12}

For $k = 6n+1$, I found evaluations of the form
\begin{equation}
\zeta(k)=E_{k,0}\sqrt3\pi^k+\sum_{d\in\{2,3,6,12\}}E_{k,d}S_k(q_3^d),\quad q_3=\exp(-\pi\sqrt3),\label{Ek}\end{equation}
with results up to $\zeta(127)$.
The first member of this family is the evaluation
\begin{equation}\zeta(7)=\tfrac{1}{1043}\left(
\tfrac{3257\sqrt3\pi^7}{16200}-2215S_7(q_3^2)-129S_7(q_3^3)
+\tfrac{129^2}{64}S_7(q_3^6)-\tfrac{129}{64}S_7(q_3^{12})\right)
\label{E7}\end{equation}
with merit $a=12$ and hence faster convergence than~(\ref{B7}), with lesser merit $a=8$.

The coefficients of the Lambert series in~(\ref{Ek}) satisfy
\begin{gather}
E_{6n+1,2}+2=E_{6n+1,3}=2^{6n}E_{6n+1,12}=-\frac{f_n}{e_n},\label{Ed}\\
E_{6n+1,6}=-f_n E_{6n+1,12},\quad f_n=2^{6n+1}+1,
\label{Ed1}\end{gather}
with $e_n={\tt 1043,\, 792701,\, 580878431,\, 423627261785,\, 308835631574603,\ldots}$\\
This integer sequence obeys a recurrence relation whose characteristic polynomial is
$(x-3^6)(x+3^3)(x-2^6)$, with $7|e_n$ and $e_0=0$, from which I inferred that
\begin{equation}e_n=\frac{3^{6n+1}-(-3)^{3n}}{2}-2^{6n}.\label{en}\end{equation}
Using the coefficient of $\sqrt3\pi^k$ in~(\ref{Ek}), I define the rational sequence
\begin{equation}E_n=\frac{(6n)!(6n+3)(6n+4)}{2^{6n-4}}e_n E_{6n+1,0}\label{Ei}\end{equation}
with $E_n={\tt 3257,\, 212373152/3,\, 325348272142978/15,\ldots}$\\
{\bf Conjecture 5:} The denominator of $E_n$ is divisible by no prime greater than $2n-1$
and by no prime congruent to 1 modulo 6.\\
{\bf Comment 5:} For $\zeta(235)$, I found that $47\cdot59\cdot75E_{39}/153734$
is a 385-digit prime.
For $\zeta(2881)$, the numerator of $E_{480}$ is divisible by a 7794-digit probable prime.

\subsection{Families with Lambert series of fractional merit}

I found 15 families that involve Lambert series $S_k(q_s^d)$ with $q_s=e^{-\pi\sqrt{s}}$,
fractional $s\in\{\frac35,\,\frac53,\,\frac13,\,\frac12,\,\frac15,\,\frac{1}{11},\frac75,\,\frac57,\,\frac17,\,\frac{1}{23},\,\frac{1}{25}\}$
and $d\sqrt{s}>2$.
They result from complex multiplication with discriminants $D\in\{-3,-4,-7,-8,-11,-15,-20,-23,-35\}$.

For $D=-15$, there are two independent families. The first begins with
\begin{align}
\zeta(3)&=\tfrac{10\sqrt{15}}{999}\pi^3
-\tfrac{3}{37}S_3(q_{5/3}^2)
+\tfrac{80}{37}S_3(q_{15})
-\tfrac{171}{37}S_3(q_{15}^2)
+\tfrac{20}{37}S_3(q_{15}^4)\label{F3}\\
\zeta(5)&=\tfrac{1327\sqrt{15}}{1516725}\pi^5
-\tfrac{63}{535}S_5(q_{5/3}^2)
+\tfrac{208}{107}S_5(q_{15})
-\tfrac{2112}{535}S_5(q_{15}^2)
+\tfrac{13}{107}S_5(q_{15}^4)\label{F5}\end{align}
with terms of merit $a\ge\frac{20}{3}$. In the denominators there is an integer sequence
$r_n={\tt 37,\, 535,\, 8047,\, 130495,\, 2103727,\, 33561055,\, 536581327,\dots}$ with solution
\begin{equation}r_n=
2^{4n+1}-\Re\left(\rho^{2n}+(-3\rho)^n\right),\quad \rho=\tfrac12(1+\sqrt{-15}).\label{rn}\end{equation}
A better family for $D=-15$ has terms of merit $a\ge\frac{48}{5}$ and begins with
\begin{equation}\zeta(3)=\tfrac{16351}{330^2\sqrt{15}}\pi^3
+\tfrac{2}{121}S_3(q_{3/5}^4)
-\tfrac{6}{605}S_3(q_{5/3}^4)
+\tfrac{342}{605}(4S_3(q_{15})+S_3(q_{15}^4))
-\tfrac{2924}{605}S_3(q_{15}^2).\label{Fa3}\end{equation}
Studying denominators of the rational coefficients, I encountered an integer sequence
$\widetilde{r}_n={\tt 605,\, 32957,\, 2107325,\, 134451197,\, 8588427965,\dots}$ with an intricate solution
\begin{gather}\widetilde{r}_n=2^{2n+1}g_n-b_n\left(5^n+(-3)^n+b_nf_n\right)\label{rnt}\\
b_n=2\Re(\rho^n),\quad f_n=2^{2n+1}+1,\quad g_n=(f_n-2^n)(f_n+2^n).\label{bfg}\end{gather}
For $k=2n+1$, the general form is
\begin{gather}\tfrac12\widetilde{r}_n\zeta(k)=
a_n\sqrt{15}\pi^k+b_n\left(5^nS_k(q_{3/5}^4)+(-3)^nS_k(q_{5/3}^4)\right)
+\sum_{d\in\{1,2,4\}}c_{n,d}S_k(q_{15}^d)\label{gen15}\\
c_{n,1}=2^{2n}c_{n,4}=f_n(g_n-b_n^2),\quad2^{2n}c_{n,2}=(b_nf_n-g_n)(b_nf_n+g_n)\label{con15} \end{gather}
with rational $a_n$. I found that $a_{1035}$ is divisible by a 6181-digit probable prime.

For $D=-3$, there is family whose terms have merit $a\ge\frac{16}{3}$, beginning with
\begin{align}
\zeta(3)&=\tfrac{133\sqrt3}{5940}\pi^3
-\tfrac{6}{11}S_3(q_{1/3}^4)
-\tfrac{18}{11}S_3(q_3^2)
+\tfrac{2}{11}S_3(q_3^4)\label{H3}\\
\zeta(5)&=\tfrac{17\sqrt3}{8694}\pi^5
+\tfrac{18}{23}S_5(q_{1/3}^4)
-\tfrac{66}{23}S_5(q_3^2)
+\tfrac{2}{23}S_5(q_3^4)\label{H5}\\
\zeta(7)&=\tfrac{10309\sqrt3}{53468100}\pi^7
-\tfrac{2322}{6601}S_7(q_{1/3}^4)
-\tfrac{10966}{6601}S_7(q_3^2)
+\tfrac{86}{6601}S_7(q_3^4).\label{H7}\end{align}
For $k=2n+1$, with $3\nmid n$, the denominator of the final coefficient is $2^{2n+1}-(-3)^n$.
For $k=6n+1$, the denominator sequence $s_n={\tt 6601,\, 20377357,\, 95065741729,\ldots}$
obeys a fourth order recurrence relation, with $7|s_n$ and $s_0=0$, solved by
\begin{equation}s_n=\tfrac13\left((2^{6n+2}-1)2^{6n}-(2^{6n+1}+1)(-3)^{3n}\right).\label{sn}\end{equation}

For $D=-8$, there is family whose terms have merit $a\ge\frac{9}{2}$, beginning with
\begin{align}
\zeta(3)&=\tfrac{181\sqrt2}{660}\pi^3
-\tfrac{16}{55}S_3(q_{1/2}^3)
-\tfrac{4}{5}S_3(q_{2}^2)
-\tfrac{52}{55}S_3(q_{2}^4)
+\tfrac{2}{55}S_3(q_{2}^{12})\label{G3}\\
\zeta(5)&=\tfrac{5171\sqrt2}{2158380}\pi^5
+\tfrac{128}{571}S_5(q_{1/2}^3)
-\tfrac{812}{571}S_5(q_{2}^2)
-\tfrac{460}{571}S_5(q_{2}^4)
+\tfrac{2}{571}S_5(q_{2}^{12})
\label{G5}.\end{align}
Studying denominators of the rational coefficients, I encountered an integer sequence
$t_n={\tt 55,\, 571,\, -3557,\, -1277,\,1212475,\, -4016045,\, -82965653,\dots}$ with solution
\begin{equation}t_n=
3^{2n+1}-(-2)^{3n}-(2^{2n+2}+(-2)^{n+1})\Re\left((1+\sqrt{-2})^{2n}\right).
\label{tn}\end{equation}

For $D=-20$, there is family whose terms have merit $a\ge5$, beginning with
\begin{align}
\zeta(3)&=\tfrac{13\sqrt5}{750}\pi^3
+\tfrac{2}{5}S_3(q_5)
-\tfrac{1}{4}S_3(q_{1/5}^6)
-\tfrac{23}{10}S_3(q_5^2)
+\tfrac{1}{10}S_3(q_5^4)
+\tfrac{1}{20}S_3(q_5^6)
\label{I3}\\
\zeta(5)&=\tfrac{809\sqrt5}{533925}\pi^5
+\tfrac{8}{113}S_5(q_5)
+\tfrac{25}{113}S_5(q_{1/5}^6)
-\tfrac{521}{226}S_5(q_5^2)
+\tfrac{1}{226}S_5(q_5^4)
+\tfrac{1}{113}S_5(q_5^6).
\label{I5}\end{align}
Studying denominators of the rational coefficients, I encountered an integer sequence
$u_n={\tt 20,\, 113,\, 980,\, 10793,\, 86360,\, 774053,\, 7471820,\, 62858993,\ldots}$ with solution
\begin{equation}u_n=\tfrac12(3^{2n+1}-(-5)^n)-\Re\left((2\sigma)^n\right),\quad \sigma=(-2+\sqrt{-5}).\label{un}\end{equation}
A better family for $D=-20$ has terms of merit $a\ge\frac{36}{5}$ and begins with
\begin{equation}\zeta(3)=\tfrac{18491\sqrt5}{1066500}\pi^3
-\tfrac{10}{79}S_3(q_{1/5}^6)
-\tfrac{2}{79}S_3(q_{1/5}^{12})
+\tfrac{18}{395}S_3(q_5^{4/3})
+\tfrac{2}{395}\sum_{d\in\{2, 3, 6, 12\}}\widetilde{g}_dS_3(q_5^d)
\label{Ia3}\end{equation}
with $\widetilde{g}_d={\tt -384,\, 4,\, 5,\, 1}$ for $d={\tt 2,3,6,12}$. The denominators yield an integer sequence
$\widetilde{u}_n={\tt 790,\, 107111,\, 20593510,\, 3025683569,\, 511718343910,\ldots}$ with a 12th 
order recurrence relation, solved by
\begin{align}\widetilde{u}_n&=2^{n-1}\left(3^{4n+2}-3^{2n+1}(-5)^n+3^{2n}-(-5)^n\right)\nonumber\\
&{}-\Re\left(\sigma^n(2^{2n+1}(-5)^n-3^{2n}+2^{2n}+(2\sigma)^n)\right).\label{unt}\end{align}

For $D=-11$, there is a family whose terms have merit $a\ge\frac{44}{9}$, beginning with
\begin{equation}\zeta(3)=\tfrac{623}{4860\sqrt{11}}\pi^3
+\tfrac{22}{45}S_3(q_{1/11}^{12})+\tfrac{2}{45}\sum_{d|36,\,d>1}f_dS_3(q_{11}^{d/3})\label{F11}\end{equation}
with integers $f_d={\tt 81,\, 132,\, -9,\, -297,\, -4,\, 33,\, 9,\, -1}$, for $d={\tt 2, 3, 4, 6, 9, 12, 18, 36}$.
Studying denominators of the rational coefficients, I encountered an integer sequence
$v_n={\tt 1,\, 29,\, 993,\, 38293,\, 1329305,\, 48535917,\, 1740077137,\ldots}$ with solution
\begin{equation}v_n=\tfrac{1}{45}\left(36^n+(-11)^n
-2\Re\left((1+\sqrt{-11})^{2n}\right)\right).\label{vn}\end{equation}
An inferior family has 21 Lambert series of merit $a\ge\frac{400}{99}$ and begins with
\begin{gather}131\,\zeta(3)=\tfrac{13\cdot3511}{2700\sqrt{11}}\pi^3
-99S_3(q_{1/11}^{20/3})-635\left(4S_3(q_{11})-9S_3(q_{11}^2)+S_3(q_{11}^4)\right)\nonumber\\{}
+\sum_{d|12,\,d>1}11\alpha_dS_3(q_{1/11}^{5d})
+\sum_{d\in\{4,6,12\}}25\beta_dS_3(q_{11}^{d/5})
+\sum_{d|36}\gamma_dS_3(q_{11}^{5d/3})
\label{alt11}\end{gather}
with integers $\alpha_d={\tt -3^223,-4,23,9,-1}$ for $d=2,3,4,6,12$,
$\beta_d={\tt -23,-9,1}$ for $d=4,6,12$,
$\gamma_d={\tt 36,-81,-2^223,9,3^223,4,-23,-9,1}$ for $d=1,2,3,4,6,9,12,18,36$.
Studying denominators of the rational coefficients, I encountered an integer sequence
$\widetilde{v}_n={\tt 262,\, 449966,\, 461557057,\, 195425426366,\ldots}$ solved for $k=2n+1$ by
\begin{align}\widetilde{v}_n&=(2g_n-1)44^n+99^n-(g_n+2^k)(-25)^n+2^k(-9)^n+b_n(c_n-5^k)(-4)^n\label{vtn}\\
g_n&=3^k+b_n+1,\; b_n=2\Re(\mu^{2n}),\; c_n=2\Re((1+\mu)^{2n}),\; \mu=\tfrac12(1+\sqrt{-11}).
\label{bcgm}\end{align}

For $D=-35$, there is a family whose terms have merit $a\ge\frac{144}{35}$, beginning with
\begin{align}\zeta(3)&=\tfrac{2969}{12960\sqrt{35}}\pi^3
+\tfrac{35}{24}S_3(q_{1/35}^{12})
-\tfrac{7}{24}S_3(q_{5/7}^{4})\nonumber\\&{}
-\tfrac{15}{8}S_3(q_{7/5}^{2})
+\tfrac{5}{24}S_3(q_{7/5}^{4})
+\tfrac{1}{24}\sum_{d|36,d>1}g_dS_3(q_{35}^{d/3})\label{got35}\end{align}
with integers $g_d={\tt 81,\, 112,\, -9,\, -252,\, -4,\, 28,\, 9,\, -1}$, for $d={\tt 2, 3, 4, 6, 9, 12, 18, 36}$.
Studying denominators of the rational coefficients, I encountered an integer sequence
$w_n={\tt 16,\, 5096,\, -140379,\, -12384944,\, 774803876,\, 34673015481,\ldots}$ with solution
\begin{equation}w_n=\tfrac13\left(36^n+(-35)^n-2\Re\left(2^{2n+1}(-5\omega)^n-(7\omega)^n\right)\right),
\quad\omega=\tfrac12(1+\sqrt{-35}).\label{wn}\end{equation}
A better family for $D=-35$ has terms of merit $a\ge\frac{28}{5}$ and begins with
\begin{gather}\zeta(3)=\tfrac{4573}{19920\sqrt{35}}\pi^3-\tfrac{315}{166}S_3(q_{7/5}^2)
+\tfrac{1}{664}(4S_3(q_{35})-9S_3(q_{35}^2)+S_3(q_{35}^4))\nonumber\\{}
+\tfrac{1}{664}(-4F_3+28F_4+9F_6-F_{12}),\quad
F_d=5S_3(q_{7/5}^d)-7S_3(q_{5/7}^d).\label{got35a}\end{gather}
Studying denominators of the rational coefficients, I encountered an integer sequence
$\widetilde{w}_n={\tt 332,\, 45583,\, 8942693,\, 1562702887,\, 284184726617,\ldots}$ with solution
\begin{equation}\widetilde{w}_n=\tfrac14\left(6(180)^n+20^n-(3^{2n+1}+2^{2n}+1)(-7)^n\right)
-\tfrac12\Re((-4\omega)^n).\label{wnt}\end{equation}

For $D=-7$ there is a family with terms of merit $a\ge\frac{36}{7}$, inferior to Family A.
At $k=3$, the result with lower merit is expressed most simply as
\begin{equation}73\,\zeta(3)=\tfrac{15}{2\sqrt7}\pi^3
-{175}S_3(q_{1/7}^6)
+{21}S_3(q_{1/7}^{12})
-{27}S_3(q_7^{4/3})
+\sum_{d|12}\widetilde{a}_dS_3(q_7^d)\label{alt7}\end{equation}
with integers $\widetilde{a}_d={\tt -12,\, 40,\, -12,\, -3,\, 25,\, -3}$ for $d={\tt 1, 2, 3, 4, 6, 12}$. Identity~(\ref{alt7}) may be combined 
with a rational multiple of~(\ref{A3}), which contains $S_3(q_7)$.

For $D=-23$, there is a family whose terms have merit $a\ge\frac{144}{23}$, beginning with
\begin{equation}\zeta(3)=\tfrac{93}{500\sqrt{23}}\pi^3-\tfrac{23}{25}S_3(q_{1/23}^{12})
+\tfrac{1}{25}\sum_{d|36,d>1}h_dS_3(q_{23}^{d/3})\label{got23}\end{equation}
with integers $h_d={\tt -72,\, -60,\, 9, 114,\, 4,\, -15,\, -8,\, 1}$, for $d={\tt 2, 3, 4, 6, 9, 12, 18, 36}$.
Studying denominators of the rational coefficients, I encountered an integer sequence
$x_n={\tt 25,\, 85132,\, -50570780,\, 66436773424,\, -39027031565300,\ldots}$ with solution
\begin{gather}x_n=\Im\left(\frac{\alpha(\alpha+2\beta-2^{2n+1}\gamma_++\gamma_-)-\beta\gamma_+}
{2^{2n+2}\sqrt{23}}\right),\\
\quad\alpha=(1+\sqrt{-23})^{2n},\quad\beta=(5-\sqrt{-23})^{2n},
\quad\gamma_{\pm}=36^n\pm(-23)^n.\label{xn}\end{gather}

For $k=4n+1$ and $D=-4$, a family with terms of merit $a\ge\frac{144}{25}$ begins with
\begin{align}
\tfrac{2639795}{74}\zeta(5)&=\tfrac{12075457}{99900}\pi^5
+5^4S_5(q_1^{12/5})-\tfrac{3050976}{37}S_5(q_1^4)\nonumber\\
&{}+2^5389S_5(q_1^6)+\sum_{d\in\{2, 4, 6, 9, 18, 36\}}\widetilde{C}_dS_5(q_1^{5d/3})\label{alt4}\end{align}
with $\widetilde{C}_d={\tt -3^437,\, 3^4,\, 2^461,\, 2^4,\, -37,\, 1}$
for $d={\tt 2, 4, 6, 9, 18, 36}$. The integer sequence
$\widetilde{c}_n={\tt 2639795,\, 598022567345,\, 123944841925119215,\ldots}$ for the denominators
in the result for $\zeta(4n+1)$
obeys a recurrence relation of order 16, with $5|\widetilde{c}_n$ and $\widetilde{c}_0=0$, solved for $k=4n+1$ by
\begin{align}\widetilde{c}_n&=4^n(15^k-10^k-5^{4n}-2^{4n})
-4(-12^4)^n+7(-10^4)^n+(-5^4)^n-2(-4^4)^n\nonumber\\
&{}-(2^k-1)18^{2n}-(2^k-(-4)^n+1)((6+8i)^{2n}+(6-8i)^{2n}).\label{cnt}\end{align}

For $k=6n+1$ and $D=-3$, a family with terms of merit $a\ge\frac{25}{3}$ begins with
\begin{equation}970\,\zeta(7)=\tfrac{3029\sqrt3}{16200}\pi^7-\tfrac{29\cdot449}{7}S_7(q_3^2)
+\tfrac{43}{2^77}\sum_{d|12}G_dS_7(q_{1/3}^{5d})\label{G12}\end{equation}
with familiar integers $G_d={\tt 12^3,\, -3^3129,\, -64,\, 27,\, 129,\, -1}$ for $d={\tt 1, 2, 3, 4, 6, 12}$.
Studying denominators of the rational coefficients, I encountered an integer sequence
$y_n={\tt 970,\, 14460979,\, 227188159336,\, 3547683233832985,\ldots}$ with solution
\begin{equation}y_n=\tfrac{1}{84}\left(5^{6n+1}-2^{6n+1}-(2^{6n+1}+1)(-3)^{3n}\right).\label{yn}\end{equation}
A better family has terms of merit $a\ge\frac{49}{3}$ and begins with
\begin{equation}44071\,\zeta(7)=\tfrac{5734\sqrt3}{675}\pi^7+\tfrac{1}{2^7}
\sum_{d_1\in\{1,2,4\}}\sum_{d_2\in\{7,9,21\}}
H_{d_1d_2}S_7(q_{1/3}^{d_1d_2})\label{H84}\end{equation}
with integers $H_d$ for $d\in\{\tt 7, 9, 14, 18, 21, 28, 36, 42, 84\}$ and
$H_7/12^3=-H_{84}=443$.
Studying denominators of the rational coefficients, I encountered an integer sequence
$z_n={\tt 44071,\, 4786291978,\, 571088556271897,\ldots}$ with solution
\begin{equation}z_n=\tfrac{1}{20}\left(7^{6n+1}-3^{6n}\left(2-(-3)^{3n+1}\right)
-2\Re((2+\sqrt{-3})^{6n})\right).\label{zn}\end{equation}

\subsection{Order of merit}

In summary, there are 20 known families with $a>4$. Their order of merit is
\begin{equation}{\tt \tfrac{49}{3} ,\, 16 ,\, 12 ,\, \tfrac{48}{5} ,\, 9 ,\,\tfrac{25}{3} ,\,
 8 ,\, \tfrac{36}{5} ,\, 7,\, \tfrac{20}{3} ,\, \tfrac{144}{23} ,\, \tfrac{144}{25} ,\, \tfrac{28}{5} ,\, 
\tfrac{16}{3} ,\, \tfrac{36}{7} ,\, 5 ,\, \tfrac{44}{9} ,\, \tfrac{9}{2} ,\, \tfrac{144}{35}} ,\,\tfrac{400}{99}\label{om}\end{equation}
in which $D=-3$ appears 4 times and $D=-4$ appears 3 times.
Each discriminant $D\in\{-7,-8,-11,-15,-20,-35\}$ appears twice, while $D=-23$ appears just once.
It seems probable that there are finitely many independent families with rational merit $a>4$.
Whether there be more than these 20 remains an open question.

\subsection{Logarithms at $k=1$}

In the singular case $k=1$, Families A to E give
\begin{align}
S_1(q_7)-2S_1(q_7^2)+S_1(q_7^4)&=-\tfrac{\sqrt7}{24}\pi+\tfrac12\log(2)\\
2S_1(q_2^2)-S_1(q_2^3)-S_1(q_2^6)&=\tfrac{5\sqrt2}{24}\pi+\tfrac14\log(2)-\log(3)\\
S_1(q_1^3)+2S_1(q_1^4)-4S_1(q_1^6)+S_1(q_1^{12})&=\tfrac{1}{24}\pi+\log(2)-\tfrac34\log(3)\\
2S_1(q_1^4)+S_1(q_1^5)-4S_1(q_1^{10})+S_1(q_1^{20})&=\tfrac{7}{24}\pi+\log(2)-\log(5)\\
S_1(q_3^2)+S_1(q_3^3)-3S_1(q_3^6)+S_1(q_3^{12})&=\tfrac{\sqrt3}{24}\pi+\tfrac13\log(2)-\tfrac5{12}\log(3).
\end{align}
Similarly, for the fractional families, the logarithms of 7, 11 and 23 also appear.
The fastest convergence occurs in the identity
\begin{equation}\tfrac{5\sqrt{3}}{36}\pi+\tfrac{13}{12}\log(3)-\log(7)=L_7-3L_9+L_{21},\quad
L_d=S_1(q_{1/3}^d)-2S_1(q_{1/3}^{2d})+S_1(q_{1/3}^{4d})\label{fast}\end{equation}
with 9 Lambert series, each no bigger than $S_1(q_{1/3}^7)<3.062\times10^{-6}$.

\section{Computational efficiency}

The evaluation~\cite{P,V}
\begin{equation}\zeta(7)=\frac{19\pi^7}{56700}
-2\sum_{n=1}^\infty\frac{1}{n^7(e^{2\pi n}-1)}\label{rz7}\end{equation}
has merit $a=4$, inferior to~(\ref{E7}) in Family E, with merit $a=12$. Yet 4 Lambert series
are needed for latter. I estimate the cost of Family E as
\begin{equation}T_E=\frac{2}{\sqrt{3}}\left(\frac12+\frac13+\frac16+\frac1{12}\right)\approx1.251
\label{TE}\end{equation}
relative to the cost of~(\ref{rz7}). The relative costs of Families A and B are
\begin{equation}
T_A=\frac{2}{\sqrt{7}}\left(1+\frac12+\frac14\right)\approx1.323,\quad
T_B=\frac{2}{\sqrt{2}}\left(\frac12+\frac13+\frac16\right)\approx1.414.\label{TAB}\end{equation}
Hence I conclude that the family beginning with Ramanujan's discovery~(\ref{rz3}) is preferable in the
case of $\zeta(4n-1)$.

In the case of $\zeta(4n+1)$, the cost of the family beginning with~(\ref{rz5}) is $1+\frac12=1.5$,
since it may be written
in terms of $S_5(q_1^2)$ and $S_5(q_1^4)$, with $q_1=e^{-\pi}$. In comparison,
the costs of Families C and D are
\begin{equation}
T_C=2\left(\frac13+\frac14+\frac16+\frac1{12}\right)\approx1.667,\quad
T_D=2\left(\frac14+\frac15+\frac1{10}+\frac1{20}\right)=1.2.\label{TCD}\end{equation}
Hence I conclude that Family D is preferable in the case of $\zeta(4n+1)$.

\subsection{Comparison with other methods}

For $\zeta(3)$, the hypergeometric method of Amdeberhan and Zeilberger~\cite{AZ} is the fastest so far exploited.
Its convenience comes from the fact that the summands are rational.
For $\zeta(5)$, there is a fast method with rational summands given by Broadhurst in Eq.~72 of~\cite{DJB1},
which consists of 22 terms of the type found by Bailey, Borwein and Plouffe (BBP), in their 4-term formula~\cite{BBP}
\begin{equation}
\pi=\sum_{k=0}^\infty\frac{1}{2^{4k}}\left(
\frac{4}{8k+1}-\frac{2}{8k+4}-\frac{1}{8k+5}-\frac{1}{8k+6}
\right).\label{BBP}\end{equation}
In the case of $\zeta(5)$, there are 8 terms with $2^{-4k}$ in their summands,
8 terms with $2^{-12k}$ and 6 terms with $2^{-20k}$. These inverse powers of 2 are
well suited to numerical computation, with 20 billion decimal digits of $\zeta(5)$ obtainable
in 25 hours on a machine with 14 cores~\cite{Y}. This record was eclipsed in 2023 using a hypergeometric formula,
discovered by Y.\ Zhao
and proved by Kam Cheong Au~\cite{Au}, which yielded 200 billion decimal digits of $\zeta(5)$ in 30 hours
on a machine with 128 cores~\cite{Y}.

Using {\tt PSLQ}~\cite{BB}, Bailey has compiled a compendium~\cite{DHB} of BBP-type identities,
including some in base 3 found in~\cite{DJB2}. None of these goes beyond polylogarithmic weight 5.
In~\cite{DJB1}, it was shown how to extend base-2 BBP identities up to weight 11,
with $\zeta(7)$, $\zeta(9)$ and $\zeta(11)$ determined by increasingly refined
combinations of binary BBP terms with products of $\pi^2$ and $\log(2)$. This base-2 polylogarithmic
ladder terminates with $\zeta(11)$ in Eq.~83 of~\cite{DJB1}. Since $\pi$ and $\log(2)$
have very fast numerical evaluations, it is probable that by diligent binary programming
of identities in~\cite{DJB1} one might evaluate odd zeta values up to $\zeta(11)$
by BBP methods that are asymptotically faster than those for the families in this article.

With a base in an algebraic number field, it is possible to reach
weight 17, using polylogarithms ${\rm Li}_{17}(\alpha^{-k})$,
with $k|630$, and products of $\pi^2$ and $\log(\alpha)$,
where $\alpha\approx1.17628$ is smallest known Salem number~\cite{BBL}.
However, ladders with irrational bases such as $\alpha$ or the golden ratio $\frac12(1+\sqrt5)$
are uncompetitive with the methods of the families identified in the current work,
where convergence is governed by a much larger transcendental constant $e^{\pi\sqrt{a}}>4071$,
with integer $a\ge7$.

I suggest that for $k=4n+1>9$ the most efficient known method for evaluating
$\zeta(k)$ is provided by Family D in~(\ref{Dk}). For $k=4n-1>11$, I suggest that the best method
is provided by the Ramanujan family~\cite{Henri,V}
\begin{equation}\zeta(4n-1)=-2S_{4n-1}(q_1^2)
-\frac{(2\pi)^{4n-1}}{2}\sum_{m=0}^{2n}(-1)^m\frac{{\cal B}_{2m}}{(2m)!}\frac{{\cal B}_{4n-2m}}{(4n-2m)!}\label{Rk}\end{equation}
where ${\cal B}_k$ is the $k$th Bernoulli number and $q_1=e^{-\pi}$. It converges faster than the method of~\cite{sumalt}
for the alternating sum $(2^{1-k}-1)\zeta(k)=\sum_{n>0}(-1)^nn^{-k}$.

\section{Resolving the families}

Consider the periodic function $P_k(z)=S_k(e^{2\pi iz})=P_k(z+1)$ with complex $z$ in the upper half-plane
$\Im z>0$. For odd $k>1$ the transformation $z\to-1/z$ gives the quasi-modular relation
\begin{equation}z^{k-1}P_k(-1/z)=P_k(z)-R_k(z)\label{qm}\end{equation}
with a correction term~\cite{Gun,V}
\begin{equation}R_k(z)=\frac{z^{k-1}-1}{2}\zeta(k)+\frac{(2\pi i)^k}{2z}\sum_{m=0}^{(k+1)/2}
z^{2m}\frac{{\cal B}_{2m}}{(2m)!}\frac{{\cal B}_{k+1-2m}}{(k+1-2m)!}\label{Rz}.\end{equation}
Moreover,~(\ref{qm}) holds at $k=1$, with
\begin{equation}R_1(z)=\frac12\log(z)+\frac{\pi i}{12}\left(z-3+\frac1{z}\right).\label{R1}\end{equation}
Since $z^{k-1}R_k(-1/z)=-R_k(z)$, the aperiodic function $M_k(z)=P_k(z)-\frac12R_k(z)$ satisfies
$z^{k-1}M_k(-1/z)=M_k(z)$. Identity~(\ref{Rk}) is then obtained from $M_{4n-1}(i)=0$ which
shows that $P_k(i)=\frac12R_k(i)$ for $k=4n-1$.

Further progress may be made by using the doubling relation
\begin{equation}P_k(z-\tfrac12)+P_k(z)=\left(2+2^{1-k}\right)P_k(2z)-2^{1-k}P_k(4z)\label{dr}\end{equation}
which holds for all real $k$. To prove this, expand the left hand side in polylogarithms to obtain
the even function $\sum_{n=1}^\infty{\rm Li}_k((-q)^n)+\sum_{n=1}^\infty{\rm Li}_k(q^n)$, with $q=e^{2\pi iz}$.
Then the right hand side results from separating cases with even and odd $n$ and using the doubling
relation ${\rm Li}_k(-x)+{\rm Li}_k(x)=2^{1-k}{\rm Li}_k(x^2)$ for polylogarithms~\cite{Lewin}.

For $k=4n+1$, one obtains a family with merit $a=4$ beginning with $\zeta(5)$ in~(\ref{rz5}) by setting $z=\frac12i$ in~(\ref{dr}).
Then $P_k(i)$ and $P_k(2i)$ appear on the right hand side. On the left,~(\ref{qm}) relates $P_k(\frac12(i-1))$
to $P_k(i+1)=P_k(i)$ and $P_k(\frac12i)$ to $P_k(2i)$. Hence $P_k(i)$ is related to $P_k(2i)$ for $k=4n+1$.

\subsection{Resolution of Family A}

For Family A in Section~2.1, I use the identity
\begin{gather}P_k(z-\tfrac14)+P_k(z+\tfrac14)=\nonumber\\
\left((2+2^{1-k})^2-2^{1-k}\right)P_k(4z)
{}-(2+2^{1-k})\left(P_k(2z)+2^{1-k}P_k(8z)\right)\label{z4}\end{gather}
obtained by two applications of~(\ref{dr}). Set $z=\frac14\sqrt{-7}$ in~(\ref{z4}). Then $S_k({q_7^d})$ with $d|4$
appear on the right hand side. On the left, the quasi-modular relation (\ref{qm}) transforms
$P_k(\frac14(\sqrt{-7}\pm1))$ to $P_k(\frac12(\sqrt{-7}\mp1))$ and ~(\ref{dr}) gives the same targets.

Cleaning up, I proved the empirical determinations of Section~2.1 and obtained
\begin{gather}A_n=\frac{(2n+3)2^{2n+3}}{\sqrt7}\Im {\cal H}_n\left(\frac{1+\sqrt{-7}}{4}\right)\label{Ais}\\
{\cal H}_n(z)=(-1)^{n+1}\sum_{m=0}^{n+1}z^{2m-1}\binom{2n+2}{2m}{\cal B}_{2m}{\cal B}_{2n+2-2m}\label{Hn}\end{gather}
for the rational sequence in Conjecture~1, now checked up to $n=1501$ for $\zeta(3003)$.

\subsection{Resolution of Family B}

To resolve Family B in~Section 2.2, I use the trebling relation
\begin{equation}P_k(z-\tfrac13)+P_k(z)+P_k(z+\tfrac13)=\left(3+3^{1-k}\right)P_k(3z)-3^{1-k}P_k(9z).\label{tr}\end{equation}
Set $z=\frac13\sqrt{-2}$ in~(\ref{tr}). Then $S_k(q_2^2)$ and $S_k(q_2^6)$ appear on the right hand side.
On the left, transform $P_k(z)$ to $S_k(q_2^3)$. Transform $P_k(\frac13(\sqrt{-2}\pm1))$ to $P_k(\sqrt{-2}\mp1)=S_k(q_2^2)$.

Cleaning up, I proved the empirical determinations of Section~2.2 and obtained
\begin{equation}B_n=\frac{4(2n+3)3^{2n}}{\sqrt2}\Im\left(
2{\cal H}_n\left(\frac{1+\sqrt{-2}}{3}\right)+{\cal H}_n\left(\frac{\sqrt{-2}}{3}\right)\right)\label{Bis}\end{equation}
for the rational sequence in Conjecture~2, now checked up to $n=1501$.

\subsection{Resolution of Family C}

To resolve Family C in Section 2.3, with $k=4n+1$, set $z=\frac23i$ in the trebling relation~(\ref{tr}).
Then $S_k(q_1^4)$ and $S_k(q_1^{12})$ appear on the right hand side.
On the left, transform $P_k(\frac23i)$ to $P_k(\frac32i)=S_k(q_1^3)$. It remains to consider
the real combination $P_k(\frac13(2i-1))+P_k(\frac13(2i+1))=P_k(\frac23(i+1))+P_k(\frac23(i-1))$,
which~(\ref{qm}) converts to the real combination $P_k(\frac34(i-1))+P_k(\frac34(i+1))$, since $(i-1)^{k-1}=(i+1)^{k-1}$.
Translate this to $P_k(\frac14(3i+1))+P_k(\frac14(3i-1))$ and use~(\ref{z4}) to obtain
$S_k(q_1^d)$ with $d\in\{3,6,12\}$.

Cleaning up, I proved the empirical determinations of Section~2.3 and obtained
\begin{equation}C_n=\frac{(4n+3)3^{4n}}{4n+1}\Im\left(
2{\cal H}_{2n}\left(\frac{2+2i}{3}\right)+{\cal H}_{2n}\left(\frac{2i}{3}\right)\right)\label{Cis}\end{equation}
for the rational sequence in Conjecture~3, now checked up to $n=750$ for $\zeta(3001)$.

\subsection{Resolution of Family D}

To resolve Family D in Section 2.4, with $k=4n+1$, set $z=\frac25i$ in the identity
\begin{equation}\sum_{r=-2}^2P_k\left(z+\tfrac{r}5\right)=
\left(5+5^{1-k}\right)P_k(5z)-5^{1-k}P_k(25z).\label{qr}\end{equation}
Then $S_k(q_1^4)$ and $S_k(q_1^{20})$ appear on the right hand side.
On the left, transform $P_k(\frac25i)$ to $S_k(q_1^5)$. Observe that $5=(1+2i)(1-2i)$
and transform $P_k(\frac15(2i\pm1))$ to $S_k(q_1^4)$.
Transform $P_k(\frac15(2i\pm2))$ to $P_k(\frac54(i\mp1))=P_k(\frac14(5i\mp1))$ and use~(\ref{z4})
to obtain $S_k(q_1^d)$ with $d\in\{5,10,20\}$.

Cleaning up, I proved the empirical determinations of Section~2.4 and obtained
\begin{equation}D_n=\frac{(4n+3)5^{4n}}{4n+1}\Im\left(
2{\cal H}_{2n}\left(\frac{2+2i}{5}\right)+2{\cal H}_{2n}\left(\frac{1+2i}{5}\right)+{\cal H}_{2n}\left(\frac{2i}{5}\right)\right)\label{Dis}\end{equation}
for the rational sequence in Conjecture~4, now checked up to $n=750$.

\subsection{Resolution of Family E}

To resolve Family E in Section 2.5, with $k=6n+1$, set $z=\frac13v$
in~(\ref{tr}), with $v=\sqrt{-3}$. Then $S_k(q_3^2)$ and $S_k(q_3^6)$ appear on the right hand side.
On the left, transform $P_k(\frac13v)$ to $S_k(q_3^2)$. Noting that $(v+1)^{k-1}=(v-1)^{k-1}$,
transform $P_k(\frac13(v\pm1))$ to $P_k(\frac34(v\mp1))=P_k(\frac14(3v\pm1))$ and use~(\ref{z4})
to obtain $S_k(q_3^d)$ with $d\in\{3,6,12\}$.

Cleaning up, I proved the empirical determinations of Section~2.5 and obtained
\begin{equation}E_n=\frac{48(2n+1)(3n+2)3^{6n}}{(3n+1)(6n+1)\sqrt3}\Im\left(
2{\cal H}_{3n}\left(\frac{1+\sqrt{-3}}{3}\right)+{\cal H}_{3n}\left(\frac{\sqrt{-3}}{3}\right)\right)\label{Eis}\end{equation}
for the rational sequence in Conjecture~5, now checked up to $n=500$ for $\zeta(3001)$.

\subsection{Resolution of fractional families}

Here more ingenuity was needed. From 3 applications of~(\ref{dr}), I proved that
\begin{gather}P_k(z-\tfrac38)+P_k(z-\tfrac18)+P_k(z+\tfrac18)+P_k(z+\tfrac38)=\nonumber\\{}
4P_k(8z)-4(1+2^{-k}+4^{-k})(P_k(4z)-(1+2^{1-k})P_k(8z)+2^{1-k}P_k(16z))\label{z8}\end{gather}
and by diligent application of~(\ref{dr},\ref{tr}) obtained the identity
\begin{gather}
P_k(z-\tfrac16)+P_k(z+\tfrac16)=P_k(z)+2^{1-k}P_k(4z)+3^{1-k}P_k(9z)\nonumber\\
{}-(2^k+1)\left(2^{1-k}P_k(2z)+6^{1-k}P_k(18z)\right)-(3^k+1)\left(3^{1-k}P_k(3z)+6^{1-k}P_k(12z)\right)\nonumber\\
{}+6^{1-k}\left((2^k+1)(3^k+1)P_k(6z)+P_k(36z)\right)\label{z6}.\end{gather}

For $D=-15$, set $z=\frac18t$ in~(\ref{z8}) with $t=\sqrt{-15}$.
Then $S_k(q_{15}^d)$ with $d|4$ appear on the right hand side.
On the left, transform $P_k(\frac18(t\pm1))$ to $P_k(\frac12(t\mp1))$,
which doubling relates to $S_k(q_{15}^d)$ with $d|4$.
Transform $P_k(\frac18(t\pm3))$ to $P_k(\frac13(t\mp3))=S_k(q_{5/3}^2)$
with merit $a=\frac{20}{3}$. The result is the family beginning with~(\ref{F3}). 
Now set $z=\frac14t$ in~(\ref{z4}).
Then $S_k(q_{15}^d)$ with $d|4$ appear on the right hand side.
On the left, transform $P_k(\frac14(t\pm1))=P_k(\frac14(t\mp3))$ to $P_k(\frac16(t\pm3))$,
which doubling converts to $S_k(q_{5/3}^d)$ with $d|4$. Transform
$S_k(q_{5/3})$ to $S_k(q_{3/5}^4)$ with merit $a=\frac{48}{5}$. Eliminate
$S_k(q_{5/3}^2)$ using the first family.
The result is a better family beginning with~(\ref{Fa3}).

For $D=-3$, one should distinguish the easier case $k=2n+1$, with $3\nmid n$, from the harder
case $k=6n+1$. For the former, set $z=\frac12\sqrt{-3}$ in the doubling relation~(\ref{dr}).
Then there is relation between $S_k(q_3^d)$ with $d|4$
and $P_k(\lambda)$
with $\lambda=\frac12(1+\sqrt{-3})$ and $\lambda^6=1$. Hence~(\ref{qm}) gives
$(1-\lambda^{k-1})P_k(\lambda)=R_k(\lambda)\ne0$ and the easier case is now resolved
by transforming $S_k(q_3)$ to $S_k(q_{1/3}^4)$.
For the harder case, set $z=\frac18\sqrt{-3}$ in~(\ref{z8}).
Then $S_k(q_3^d)$ with $d|4$ appear on the right hand side.
On the left, transform $P_k(\frac18(\sqrt{-3}\pm1))$ to $S_k(q_3^4)$.
Transform $P_k(\frac18(\sqrt{-3}\pm3))$ to $S_k(q_{1/3}^4)$.

For $D=-8$, set $z=\frac23\sqrt{-2}$ in~(\ref{tr}).
Then $S_k(q_2^4)$ and $S_k(q_2^{12})$ appear on the right hand side. On the left, transform
$P_k(z)$ to $S_k(q_{1/2}^3)$. Transform $P_k(z\pm\frac13)=P_k(z\mp\frac23)$ to $P_k(\frac12(\sqrt{-2}\pm1))$,
which doubling relates to $S_k(q_2^d)$ with $d|4$. Transform $S_k(q_2)$ to $S_k(q_2^2)$.

For $D=-20$, set $z=\frac13\sqrt{-5}$ in~(\ref{tr}).
Then $S_k(q_5^2)$ and $S_k(q_5^6)$ appear on the right hand side. On the left, transform
$P_k(z)$ to $S_k(q_{1/5}^6)$. Transform $P_k(\frac13(\sqrt{-5}\pm1))$ to $P_k(\frac12(\sqrt{-5}\mp1))$,
which doubling relates to $S_k(q_5 ^d)$ with $d|4$. The result is the family beginning
with~(\ref{I3}). Now set $z=\frac16\sqrt{-5}$ in~(\ref{z6}). 
Then $S_k(q_5^{d/3})$ with $d|36$ appear on the
right hand side. Transform $S_k(q_5^{1/3})$ to $S_k(q_{1/5}^{12})$.
Transform $S_k(q_5^{2/3})$ to $S_k(q_{1/5}^{6})$. Eliminate $S_k(q_5)$ using the first family.
On the left, transform $P_k(\frac16(\sqrt{5}\pm1))$ to $S_k(q_5^2)$. The result is the family beginning
with~(\ref{Ia3}), with terms of merit $a\ge\frac{36}{5}$.

For $D=-11$, set $z=\frac16t$ with $t=\sqrt{-11}$ in~(\ref{z6}).
Then $S_ k(q_{11}^{d/3})$ with $d|36$ appear on the right hand side.
Transform $S_k(q_{11}^{1/3})$ to $S_k(q_{1/11}^{12})$.
On the left, transform $P_k(\frac16(t\pm1))$ to $P_k(\frac12(t\mp1))$, which doubling
relates to $S_k(q_{11}^d)$ with $d|4$. The result is the family beginning with~(\ref{F11}),
with terms of merit $a\ge\frac{44}{9}$. To resolve the inferior family beginning with~(\ref{alt11}),
with 21 Lambert series of merit $a\ge\frac{400}{99}>4$, I combined~(\ref{qr},\ref{z6}) and derived an identity with 35 terms,
exploiting the group of units modulo 30, with residues $\pm r\in\{1,7,11,13\}$.
On the left are the 8 terms $P_k(z+\frac{r}{30})$. On the right are the 27 terms $P_k(dz)$ with $d|900$.
With $z=\frac{1}{30}t$, the right hand side gives $S_k(q_{11}^{5d/3})$, $S_k(q_{11}^{d/3})$
and $S_k(q_{11}^{d/15})$, with $d|36$. Transform $S_k(q_{11}^{1/3})$ to $S_k(q_{1/11}^{12})$.
For $d=1, 2, 3, 4, 6, 9$, transform $S_k(q_{11}^{d/15})$ 
to $S_k(q_{1/11}^{60/d})$. Then all 27 terms on the right have merit $a\ge\frac{400}{99}$.
On the left, transform $P_k(\frac1{30}(t\pm1))$ to
$P_k(\frac52(t\mp1))$, which doubling relates to $S_k(q_{11}^{5d})$ with $d|4$.
Transform $P_k(\frac1{30}(t\pm7))$ to $P_k(\frac12(t\mp7))$, 
which doubling relates to $S_k(q_{11}^d)$ with $d|4$.
Transform $P_k(\frac1{30}(t\pm11))$ to $P_k(\frac5{22}(t\mp11))$, 
which doubling relates to $S_k(q_{1/11}^{5d})$ with $d|4$. 
Transform $S_k(q_{1/11}^5)$ to $S_k(q_{11}^{4/5})$.
Transform $P_k(\frac1{30}(t\pm13))$ to $P_k(\frac16(t\mp13))=P_k(\frac16(t\pm5))$.
Observe that~(\ref{qm}) relates $P_k(\frac16(t+5))$
to $P_k(\frac16(t-5))$ and that~(\ref{z6}) relates
$P_k(\frac16(t+5))+P_k(\frac16(t-5))$ 
to $S_k(q_{11}^{d/3})$ with $d|36$. Then all of the Lambert series on the left also occur on the right.
Finally, reduce these 27 terms to 21, by using the superior family to eliminate $S_k(q_{1/11}^{12})$.
With so many transformations in play, the denominator sequence~(\ref{vtn}) for the inferior family
is rather complicated, obeying a recurrence relation of order 17.
 
For $D=-35$, set $z=\frac16v$ in~(\ref{z6}) with $v=\sqrt{-35}$.
Then $S_ k(q_{35}^{d/3})$ with $d|36$ appear on the right hand side. Transform
$S_k(q_{35}^{1/3})$ to $S_k(q_{1/35}^{12})$.
On the left, transform $P_k(\frac16(v\pm1))=P_k(\frac16(v\mp5))$ to $P_k(\frac1{10}(v\pm5))$,
which doubling relates to $S_k(q_{7/5}^d)$ with $d|4$. Transform $S_k(q_{7/5})$ to
$S_k(q_{5/7}^4)$ to obtain the family with terms of merit $a\ge\frac{144}{35}$ beginning with~(\ref{got35}).
Now set $z=\frac16w$ in~(\ref{z6}) with $w=\sqrt{-7/5}$.
On the right hand side, transform $S_ k(q_{7/5}^{d/3})$ with $d\in\{1,2,3,4\}$
to obtain $S_ k(q_{5/7}^d)$ with $d\in\{3,4,6,12\}$.
On the left, transform $P_k(\frac16(w\pm1))$ to $P_k(\frac52(w\mp1))$,
which doubling relates to $S_k(q_{35}^d)$ with $d|4$.
This results in a family with terms of merit $a\ge\frac{28}{5}$, beginning with~(\ref{got35a}).

For $D=-7$, set $z=\frac{1}{12}\sqrt{-7}$ in the identity
\begin{gather}
P_k(z-\tfrac5{12})+P_k(z-\tfrac1{12})+P_k(z+\tfrac1{12})+P_k(z+\tfrac5{12})=\nonumber\\
(2^k+1)\left(2^{1-k}P_k(2z)+4^{1-k}P_k(8z)+6^{1-k}P_k(18z)+12^{1-k}P_k(72z)\right)\nonumber\\
{}-\left((2^k+1)^2-2^{k-1}\right)\left(
4^{1-k}P_k(4z)-12^{1-k}(3^k+1)P_k(12z)+12^{1-k}P_k(36z)\right)\nonumber\\
{}-6^{1-k}(2^k+1)(3^k+1)\left(P_k(6z)+2^{1-k}P_k(24z)\right)
\label{z12}\end{gather}
obtained by diligent application of~(\ref{dr},\ref{tr}).
Then $S_ k(q_{7}^{d/3})$ with $d|36$ appear on the right hand side.
Transform $S_k(q_7^{1/3})$ to $S_k(q_{1/7}^{12})$.
Transform $S_k(q_7^{2/3})$ to $S_k(q_{1/7}^6)$, with merit $a=\frac{36}{7}$.
On the left, transform $P_k(\frac{1}{12}(\sqrt{-7}\pm1))$ to 
$P_k(\frac{3}{2}(\sqrt{-7}\mp1))$, which doubling relates to $S_k(q_7^{3d})$ with $d|4$.
Transform $P_k(\frac{1}{12}(\sqrt{-7}\pm5))=P_k(\frac{1}{12}(\sqrt{-7}\mp7))$ to 
$P_k(\frac{3}{2}(\sqrt{-1/7}\pm1))$, which doubling relates to $S_k(q_{1/7}^{3d})$ with $d|4$.
Transform $S_k(q_{1/7}^{3})$ to $S_k(q_7^{4/3})$. The result is a relation with 9 terms of merit $a\ge\frac{36}{7}$
as in identity~(\ref{alt7}).

For $D=-23$, set $z=\frac16\sqrt{-23}$ in~(\ref{z6}).
Then $S_ k(q_{23}^{d/3})$ with $d|36$ appear on the right hand side. On the left, transform
$P_k(\frac16(\sqrt{-23}\pm1))$ to $P_k(\frac14(\sqrt{-23}\mp1))$. This pair occurs with complex coefficients.
Their real sum is also related to targets $S_ k(q_{23}^{d/3})$ with $d|36$ by setting $z=\frac14\sqrt{-23}$ in~(\ref{z4}).
Thus both of $P_k(\frac14(\sqrt{-23}\pm1))$ are related to these 9 targets. Now set $z=\frac1{12}\sqrt{-23}$ in~(\ref{z12}).
Then $S_ k(q_{23}^{d/3})$ with $d|36$ appear on the right hand side. On the left, transform
$P_k(z\pm\frac{1}{12})$ to $P_k(\frac12(\sqrt{-23}\mp1))$, which doubling relates to
$S_k(q_{23}^d)$ with $d|4$.
Transform $P_k(z\pm\frac{5}{12})$ to $P_k(\frac14(\sqrt{-23}\mp5))=P_k(\frac14(\sqrt{-23}\mp1))$.
It follows that there is a relation between $S_ k(q_{23}^{d/3})$ with $d|36$. At $d=1$, transform
$S_ k(q_{23}^{1/3})$ to $S_ k(q_{1/23}^{12})$, with merit $a=\frac{144}{23}$,
to obtain the family beginning with~(\ref{got23}). Sequence~(\ref{xn}) comes from a solving a pair of
simultaneous equations with complex coefficients.

For $D=-4$ and $k=4n+1$, set $z=\frac35i$ in~(\ref{qr}). Then $S_k(q_1^6)$ and $S_k(q_1^{30})$ appear on the right.
On the left, transform $P_k(z)$ to $S_k(q_1^{10/3})$. Transform $P_k(z\pm\frac15)$ to $P_k(\frac12(3i\mp1))$,
which doubling relates to $S_k(q_1^{3d})$ with $d|4$. Transform $P_k(z\pm\frac25)=P_k(z\mp\frac35)$ to 
$P_k(\frac56(i\pm1))=P_k(\frac16(5i\mp1))$, which~(\ref{z6}) relates to $S_k(q_1^{5d/3})$ with $d|36$.
Transform $S_k(q_1^{5/3})$ to $S_k(q_1^{12/5})$ with merit $a=\frac{144}{25}$. Use Family C to eliminate $S_k(q_1^3)$.
Use Family D to eliminate $S_k(q_1^5)$. The result is a family beginning with~(\ref{alt4}). 

For $D=-3$ and $k=6n+1$, set $z=\frac{1}{10}\sqrt{-3}$ in the identity
\begin{gather}\sum_{r\in\{1,3\}}\left(P_k(z-\tfrac{r}{10})+P_k(z+\tfrac{r}{10})\right)=P_k(z)+2^{1-k}P_k(4z)+5^{1-k}P_k(25z)\nonumber\\
{}-(2^k+1)\left(2^{1-k}P_k(2z)+10^{1-k}P_k(50z)\right)-(5^k+1)\left(5^{1-k}P_k(5z)+10^{1-k}P_k(20z)\right)\nonumber\\
{}+10^{1-k}\left((2^k+1)(5^k+1)P_k(10z)+P_k(100z)\right)\label{z10}.\end{gather}
Then $S_k(q_3^{d/5})$ with $d|100$ appear on the right hand side. Transform terms with $d|4$
to obtain $S_k(q_{1/3}^{5d})$ with $d|4$. At $d=5$, transform $S_k(q_3)$ to $S_k(q_{1/3}^4)$ with merit $a=\frac{16}{3}$.
Eliminate this using the family that begins with~(\ref{H7}). All the terms on the right now have merit $a\ge\frac{25}{3}$.
On the left, transform $P_k(\frac1{10}(\sqrt{-3}\pm1))$ to $P_k(\frac52(\sqrt{-3}\mp1))$, which doubling
relates to $S_k(q_3^{5d})$ with $d|4$. Transform $P_k(\frac1{10}(\sqrt{-3}\pm3))$ to $P_k(\frac56(\sqrt{-3}\mp3))$, which doubling
relates to $S_k(q_{1/3}^{5d})$ with $d|4$. The result is the family beginning with~(\ref{G12}).
To resolve a better family, with terms of merit $a\ge\frac{49}{3}$, set $z=\frac17\sqrt{-3}$ in
\begin{equation}\sum_{r=-3}^3P_k\left(z+\tfrac{r}7\right)=
\left(7+7^{1-k}\right)P_k(7z)-7^{1-k}P_k(49z)\label{sr}\end{equation}
with $S_k(q_3^2)$ and $S_k(q_3^{14})$ on the right hand side. On the left, transform $P_k(z)$ to $S_k(q_{1/3}^{14})$.
Transform $P_k(z\pm\frac17)$ to $P_k(\frac74\sqrt{-3}\pm\frac14)$ and use~(\ref{z4}) to produce
$S_k(q_3^{7d})$ with $d|4$. Transform $P_k(z\pm\frac27)$ to $P_k(\sqrt{-3}\mp2)=S_k(q_3^2)$.
Transform $P_k(z\pm\frac37)$ to $P_k(\frac{7}{12}(\sqrt{-3}\mp3))=P_k(\frac{7}{12}\sqrt{-3}\pm\frac14)$ and use~(\ref{z4}) to produce
$S_k(q_{1/3}^{7d})$ with $d|4$. Thus $S_k(q_3^2)$ is related to terms of the form
$S_k(q_3^{7d})$ and $S_k(q_{1/3}^{7d})$ with $d|4$. It is also related, by Family E in~(\ref{Ek}),
to $S_k(q_3^{3d})$ with $d|4$. By this method, I obtained relation~(\ref{H84}) between
$\zeta(7)$, $\sqrt{3}\pi^7$ and the 9 terms $S_7(q_{1/3}^{d_1d_2})$
with $d_1\in\{1,2,4\}$ and $d_2\in\{7,9,21\}$.
These 9 terms are less efficient than the 4 terms in Family E.
The computational cost of the former, relative to the latter, is $19/13$.

\subsection{Derivative relations}

Each family may be extended by taking derivatives. For odd $k>1$, let
\begin{equation}
T_k(z)=P_k(z)-\frac{2zP^\prime_k(z)}{k-1},\quad
U_k(z)=R_k(z)-\frac{2zR^\prime_k(z)}{k-1}.\label{TU}\end{equation}
Then $z^{k-1}U(-1/z)=U_k(z)$ and there is a quasi-modular relation
\begin{equation}z^{k-1}T_k(-1/z)=U_k(z)-T_k(z).\label{qmd}\end{equation}
Hence $T_k(i)=\frac12U_k(i)$ for $k=4n+1$, as recorded by Emil Grosswald,
in Theorem A of~\cite{G}, and by Henri Cohen, in Theorem 1.1 of~\cite{Henri}.
With $z/i=x$, the identity $z^2P_1^\prime(z)=z^2R_1^\prime(z)+P_1^\prime(-1/z)$ gives
\begin{equation}x^2{\cal S}(x)=\frac{1+x^2}{6}-\frac{x}{\pi}-{\cal S}\left(\frac{1}{x}\right),\quad
{\cal S}(x)=\sum_{n=1}^\infty\frac{1}{\sinh^2(n\pi x)}\label{Sx}\end{equation}
for $\Re x>0$
and hence the evaluation $\frac{1}{\pi}=\frac13-2{\cal S}(1)$ with merit $a=4$.

In Families A and B, derivative identities yield alternative methods of evaluating $\zeta(2n+1)$.
They extend Families C and D to $\zeta(4n-1)$ and Family E to $\zeta(6n-3)$. 
Similar extensions apply to the 15 families with Lambert series of fractional merit.
At $k=3$, the derivative identities for Families A to E are
\begin{align}
\zeta(3)&=\tfrac{7\sqrt7}{480}\pi^3
+\tfrac52T_3(\tfrac12\sqrt{-7})
-\tfrac{41}8T_3(\sqrt{-7})
+\tfrac58T_3(2\sqrt{-7})\\
&=\tfrac{47\sqrt2}{1710}\pi^3
-\tfrac{44}{19}T_3(\sqrt{-2})
+\tfrac{4}{19}T_3(\tfrac32\sqrt{-2})
+\tfrac{2}{19}T_3(3\sqrt{-2})\\
&=\tfrac{7}{180}\pi^3-4T_3(\tfrac32\sqrt{-1})-8T_3(2\sqrt{-1})+11T_3(3\sqrt{-1})-T_3(6\sqrt{-1})\\
&=\tfrac{221}{5700}\pi^3
-\tfrac{232}{95}T_3(2\sqrt{-1})
-\tfrac{28}{95}T_3(\tfrac52\sqrt{-1})
+\tfrac{77}{95}T_3(5\sqrt{-1})
-\tfrac{7}{95}T_3(10\sqrt{-1})\\
&=\tfrac{43\sqrt3}{1920}\pi^3
-\tfrac{25}{8}T_3(\sqrt{-3})
-\tfrac{9}{8}T_3(\tfrac32\sqrt{-3})
+\tfrac{81}{32}T_3(3\sqrt{-3})
-\tfrac{9}{32}T_3(6\sqrt{-3}).
\end{align}

Using {\tt lindep} in an extensive and systematic search,
I found 53 independent identities relating $\zeta(3)$
and an algebraic multiple of $\pi^3$ to rational multiples
of $T_3(z)$, with $z$ a rational multiple of $\sqrt{-s}$
and $s\in[1,77]$ a square-free integer.
In each of these cases there is also an identity relating 
logs of primes $p\le23$ and an algebraic multiple of $\pi$ to corresponding values of $P_1(z)$.
However, only 20 of these cases, with $s\in[1,35]$, yield $\zeta(k)$ with $k>3$,
producing the order of merit already recorded in~(\ref{om}).
The other 33 cases are isolated results that do not produce infinite families.

Using the quasi-modular transformation~(\ref{Sx}), I obtained
\begin{align}
\frac{11}{\sqrt3}\left(\frac{1}{\pi}-\frac{7}{22}\right)&=4(1+\sqrt3){\cal S}(\sqrt3)
-3(6+\sqrt3)({\cal S}(\tfrac32\sqrt3)-4{\cal S}(6\sqrt3))\nonumber\\
&{}-8(5-2\sqrt3){\cal S}(2\sqrt3)-18(1+2\sqrt3){\cal S}(3\sqrt3)\label{pi12}\end{align}
with merit $a=12$. Higher merit is achieved in
\begin{equation}\frac{11}{\pi}-\frac{7}{2}=\sum_{n=1}^{\infty}\sum_{d\in\{3,4,6,12\}}
\frac{\beta_d+\gamma_d\sqrt{15}}{\sinh^2(n\pi\sqrt{15}d/6)}
=\sum_{n=1}^{\infty}\sum_{d\in\{3,4,6,12\}}
\frac{\widetilde{\beta}_d+\widetilde{\gamma}_d\sqrt{8}}{\sinh^2(n\pi\sqrt{8}d/4)}\label{faster}\end{equation}
with $\beta_d={\tt -1,\, -24,\, 72,\, -68}$ and $\gamma_d={\tt 18,\, -8,\, -82,\, 72}$ for $d=3,4,6,12$
giving merit $a=15$, while
$\widetilde{\beta}_d={\tt 105,\, 96,\, -486,\, 264}$
and $\widetilde{\gamma}_d={\tt 39,\, 0,\, -183,\, 144}$ give $a=18$. 
Better yet are the evaluations with $a=20$ and $a=\frac{144}{7}$ in
\begin{align}\frac{1511}{\pi}&=481
-\sum_{n=1}^{\infty}\left(\frac{396+360\sqrt{5}}{\sinh^2(6n\pi/\sqrt5)}
+\sum_{d|12,\,d>1}\frac{b_d+c_d\sqrt{5}}
{\sinh^2(n\pi\sqrt{5}d/2)}\right)\label{best5}\\
\frac{111}{2\pi}&=\frac{53}{3}
-\sum_{n=1}^{\infty}\left(\frac{32+16\sqrt{7}}{\sinh^2(6n\pi/\sqrt7)}
+\sum_{d|12,\,d>1}\frac{\widetilde{b}_d+\widetilde{c}_d\sqrt{7}}
{\sinh^2(n\pi\sqrt{7}d/2)}\right)\label{best7}\end{align}
with $b_d={\tt 549,\, 1368,\, -4440,\, -2439,\, 7452}$,
$c_d={\tt 4493,\, -1080,\, -1280,\, 27,\, -2520}$,
$\widetilde{b}_d={\tt 752,\, -350,\, -1504,\, 0,\, 1176}$,
$\widetilde{c}_d={\tt 147, \, 35,\, 180,\, -630,\, 252}$
for $d=2,3,4,6,12$.
The fastest convergence occurs with merit $a=\frac{108}{5}$ in a relation
between $\frac{1}{\pi}-\frac{127}{399}$ and algebraic multiples of 8 series of the form 
${\cal S}(r\sqrt{15})$ with $r\in\{\frac35, \frac23, 1, \frac65, \frac32, 2, 3, 6\}$. In this
case the asymptotic ratio of terms $e^{-6\pi\sqrt{3/5}}<\frac{1}{2193048}$
provides at least 6 extra decimal digits of precision from each additional term.

\section{Comments and conclusions}

On 12 December 2023, Simon Plouffe kindly informed me by email of his notable
discoveries~(\ref{A5},\ref{B3},\ref{D5}), which are members of Families A, B and D.
Setting up a systematic search with {\tt lindep} in {\tt Pari/GP}, I discovered
Families C and E. Further search strongly suggested that there are only 5 independent families of relations
with integer merit $a\ge7$.
Intrigued by this apparent oligopoly, reminiscent of the mafiosi of New York in 1931,
I decided to investigate the 5 known families, up to $\zeta(127)$, seeking
explicit formul{\ae} for their rational coefficients.
This resulted in recurrence relations
for intriguing integer sequences, solved in~(\ref{an},\ref{bn},\ref{cn},\ref{dn},\ref{en}).
Then~(\ref{Rz}) gave formul{\ae}
in~(\ref{Ais},\ref{Bis},\ref{Cis},\ref{Dis},\ref{Eis}) for coefficients of $\pi^k$.
In the first instance, this was done by guesswork, based on
correspondences between the $\zeta(k)$ and $\pi^k$ terms in corrections to modularity for Lambert series.
These successful guesses for the arguments of the function ${\cal H}_n(z)$ in~(\ref{Hn})
were the key to subsequent proofs, based on
$\sum_{r\in{\mathbf Z}_n^\times}P_k(z+\frac{r}{n})$ for groups of units modulo 
$n\in\{2,3,4,5,6,7,8,10,12,30\}$.

Here I was guided by work by Spencer Bloch, Matt Kerr and Pierre Vanhove~\cite{B1,B2}
on Feynman integrals in quantum field theory, where one encounters modular forms that
can be written as infinite sums of ${\rm Li}_k(q^n)$ with $k=2$ in~\cite{B1} and $k=3$ in~\cite{B2}.
In the latter case, Detchat Samart was able to prove a conjecture~\cite{DJB3} of mine by
using the quasi-modular transformation~(\ref{qm}) and doubling relation~(\ref{dr}) at $k=3$,
in Proposition~1 of~\cite{S}. Thus prepared, I was able to construct the proofs in Section~4,
with discriminants $D\in\{-3,\,-4,\,-7,\,-8,\,-11,\,-15,\,-20,\,-23,\,-35\}$

At heart, these proofs use elementary
decompositions of divisors of 24, namely
$2=1+1$, $3=2+1$, $2^2=3+1=5-1$, $6=5+1$, $2^3=7+1=5+3$, $12=11+1=7+5$
and $24=23+1$.
As Kronecker remarked, {\em die ganzen Zahlen hat der liebe Gott gemacht, alles andere ist Menschenwerk.} 

In conclusion, I remark that the empirical results of Section 2 are now proved, with the exception of Conjectures
1 to 5, on primes in the denominators of rational numbers produced by sums over products of Bernoulli numbers,
and my conjecture that precisely 5 independent families achieve integer merit $a\ge7$.
I leave open the question as to whether more than 20 independent families have rational merit $a>4$.

\subsection*{Acknowledgements} I thank Simon Plouffe, for correspondence on Families A, B and D,
Tate Broadhurst, for a telling observation made during our joint scrutiny of Family B, and Wadim Zudilin,
for kind encouragement to discover families with terms of fractional merit.

}\raggedright


\begin{thebibliography}{99}

\bibitem{AZ}
Tewodros Amdeberhan and Doron Zeilberger,
Hypergeometric series acceleration via the WZ method,
Elec.\ J.\ Combin.\ {\bf 4} (1997) R3.

\bibitem{Au}
Kam Cheong Au,
Wilf-Zeilberger seeds and non-trivial hypergeometric identities,
arXiv:2312.14051.

\bibitem{DHB}
David H. Bailey,
A compendium of BBP-type formulas for mathematical constants,
\url{https://www.davidhbailey.com/dhbpapers/bbp-formulas.pdf} (2023).

\bibitem{BB}
David H.\ Bailey and David Broadhurst,
Parallel integer relation detection: techniques and applications,
Math.\ Comp. {\bf 70} (2001) 1719--1736.

\bibitem{BBL}
David H.\ Bailey and David Broadhurst,
A seventeenth-order polylogarithm ladder,
arXiv:math/9906134.

\bibitem{BBP}
David H.\ Bailey, Peter B.\ Borwein and Simon Plouffe,
On the rapid computation of various polylogarithmic constants,
Math.\ Comp. {\bf 66} (1997) 903--913.

\bibitem{Bruce}
Bruce Berndt,
Modular transformations and generalizations of several formul{\ae} of Ramanujan,
Rocky Mountain J.\ Math.\ {\bf 7} (1977) 147--189.

\bibitem{B1}
Spencer Bloch and Pierre Vanhove,
The elliptic dilogarithm for the sunset graph,
J.\ Number Theory {\bf 148} (2015) 328--364,
arXiv:1309.5865.

\bibitem{B2}
Spencer Bloch, Matt Kerr and Pierre Vanhove,
A Feynman integral via higher normal functions,
Compositio Math.\ {\bf 151} (2015) 2329--2375,
arXiv:1406.2664.

\bibitem{DJB1}
David Broadhurst,
Polylogarithmic ladders, hypergeometric series and the ten
millionth digits of $\zeta(3)$ and $\zeta(5)$,
arXiv:math/9803067.

\bibitem{DJB2}
David Broadhurst,
Massive 3-loop Feynman diagrams reducible to {\tt SC*}
primitives of algebras of the sixth root of unity,
Eur.\ Phys.\ J.\ {\bf C8} (1999) 311--333,
arXiv:hep-th/9803091.

\bibitem{DJB3}
David Broadhurst,
Multiple zeta values and modular forms in quantum field theory,
in {\em Computer algebra in quantum field theory}, eds.\ C.\ Schneider
and J.\ Bl\"umlein (Springer, Wien, 2013) 33--72.

\bibitem{Henri}
Henri Cohen,
High precision computation of Hardy-Littlewood constants,
\url{http://oeis.org/A221712/a221712.pdf} (1998).

\bibitem{sumalt}
Henri Cohen, Fernando Rodriguez Villegas and Don Zagier,
Convergence acceleration of alternating series,
Experiment.\ Math.\ {\bf 9} (2000) 3--12.

\bibitem{G}
Emil Grosswald,
Remarks concerning the values of the Riemann zeta function at integral, odd arguments,
J.\ Number Theory {\bf 4} (1972) 225--235.

\bibitem{Gun}
Sanoli Gun, M. Ram Murty and Purusottam Rath,
Transcendental values of certain Eichler integrals,
Bull.\ London Math.\ Soc.\ {\bf 43} (2011) 939--952.

\bibitem{Lewin}
Leonard Lewin,
{\em Polylogarithms and associated functions}, North-Holland, Amsterdam, 1981.

\bibitem{GP}
PARI Group,
{\sc Pari/GP} version 2.15.4, Univ.\ Bordeaux,\\
\url{http://pari.math.u-bordeaux.fr/} (2023).

\bibitem{P}
Simon Plouffe,
Identities inspired from the Ramanujan notebooks, first series,
arXiv:1101.4826.

\bibitem{P1}
Simon Plouffe,
Efficient formulas for $\zeta(n)$, email, 12 December 2023.

\bibitem{S}
Detchat Samart,
Feynman integrals and critical modular L-values,
Commun. Number Theory and Physics
{\bf 10} (2016) 133--156,
arXiv:1511.07947.

\bibitem{V}
Linas Vep\v{s}tas,
On Plouffe’s Ramanujan identities,
Ramanujan J.\ {\bf 7} (2012) 387--408.

\bibitem{Y}
Alex Yee,
Records set by y-cruncher,
\url{http://www.numberworld.org/y-cruncher/records/2018_10_5_zeta5.txt} (2018),\\
\url{http://www.numberworld.org/y-cruncher/records/2023_12_10_zeta5.txt} (2023).

\end{thebibliography}
\end{document}